\documentclass{amsart}
\usepackage{amssymb,amsmath}
\usepackage[mathscr]{euscript}
\input{epsf}

\newcounter{sec}

\newcounter{punct}[sec]

\def\punct{\refstepcounter{punct}{\arabic{sec}.\arabic{punct}.  }}

\newtheorem{theorem}{Theorem}[sec]
\newtheorem{proposition}[theorem]{Proposition}

\newtheorem{lemma}[theorem]{Lemma}

\newtheorem{corollary}[theorem]{Corollary}

\newtheorem{definition}[theorem]{Definition}
\newtheorem{propositionA}{Proposition A.}
\newtheorem{questionA}[propositionA]{Question A.}

\def\COUNTERS{\addtocounter{sec}{1}
              \setcounter{punct}{0}
          \setcounter{equation}{0}
          \setcounter{theorem}{0}
          }
          
          \def\sm{\smallskip}

\begin{document}

\newcommand{\supp}{\mathop {\mathrm {supp}}\nolimits}
\newcommand{\rk}{\mathop {\mathrm {rk}}\nolimits}
\newcommand{\Aut}{\mathop {\mathrm {Aut}}\nolimits}
\newcommand{\Out}{\mathop {\mathrm {Out}}\nolimits}
\renewcommand{\Re}{\mathop {\mathrm {Re}}\nolimits}

\newcommand{\GL}{\mathop {\mathrm {GL}}\nolimits}
\newcommand{\Sp}{\mathop {\mathrm {Sp}}\nolimits}
\newcommand{\SO}{\mathop {\mathrm {SO}}\nolimits}
\renewcommand{\O}{\mathop {\mathrm {O}}\nolimits}
\newcommand{\U}{\mathop {\mathrm {U}}\nolimits}
\newcommand{\SU}{\mathop {\mathrm {SU}}\nolimits}
\newcommand{\Ad}{\mathop {\mathrm {Ad}}\nolimits}
\newcommand{\Isom}{\mathop {\mathrm {Isom}}\nolimits}
\newcommand{\Abs}{\mathop {\mathrm {Abs}}\nolimits}
\newcommand{\Clop}{\mathop {\mathrm {Clop}}\nolimits}
\newcommand{\spike}{\mathop {\mathrm {spike}}\nolimits}
\renewcommand{\vert}{\mathop {\mathrm {vert}}\nolimits}
\newcommand{\Hie}{\mathop {\mathrm {Hier}}\nolimits}
\newcommand{\SL}{\mathop {\mathrm {SL}}\nolimits}
\newcommand{\PSL}{\mathop {\mathrm {PSL}}\nolimits}
\newcommand{\Ba}{\mathop {\mathrm {Ba}}\nolimits}
\newcommand{\Br}{\mathop {\mathrm {Br}}\nolimits}
\newcommand{\Diff}{\mathop {\mathrm {Diff}}\nolimits}

\def\Ss{\mathrm{S}}

\def\ov{\overline}
\def\wh{\widehat}
\def\wt{\widetilde}

\renewcommand{\rk}{\mathop {\mathrm {rk}}\nolimits}
\renewcommand{\Aut}{\mathop {\mathrm {Aut}}\nolimits}
\renewcommand{\Re}{\mathop {\mathrm {Re}}\nolimits}
\renewcommand{\Im}{\mathop {\mathrm {Im}}\nolimits}
\newcommand{\sgn}{\mathop {\mathrm {sgn}}\nolimits}

\def\bfa{\mathbf a}
\def\bfb{\mathbf b}
\def\bfc{\mathbf c}
\def\bfd{\mathbf d}
\def\bfe{\mathbf e}
\def\bff{\mathbf f}
\def\bfg{\mathbf g}
\def\bfh{\mathbf h}
\def\bfi{\mathbf i}
\def\bfj{\mathbf j}
\def\bfk{\mathbf k}
\def\bfl{\mathbf l}
\def\bfm{\mathbf m}
\def\bfn{\mathbf n}
\def\bfo{\mathbf o}
\def\bfp{\mathbf p}
\def\bfq{\mathbf q}
\def\bfr{\mathbf r}
\def\bfs{\mathbf s}
\def\bft{\mathbf t}
\def\bfu{\mathbf u}
\def\bfv{\mathbf v}
\def\bfw{\mathbf w}
\def\bfx{\mathbf x}
\def\bfy{\mathbf y}
\def\bfz{\mathbf z}

\def\bfA{\mathbf A}
\def\bfB{\mathbf B}
\def\bfC{\mathbf C}
\def\bfD{\mathbf D}
\def\bfE{\mathbf E}
\def\bfF{\mathbf F}
\def\bfG{\mathbf G}
\def\bfH{\mathbf H}
\def\bfI{\mathbf I}
\def\bfJ{\mathbf J}
\def\bfK{\mathbf K}
\def\bfL{\mathbf L}
\def\bfM{\mathbf M}
\def\bfN{\mathbf N}
\def\bfO{\mathbf O}
\def\bfP{\mathbf P}
\def\bfQ{\mathbf Q}
\def\bfR{\mathbf R}
\def\bfS{\mathbf S}
\def\bfT{\mathbf T}
\def\bfU{\mathbf U}
\def\bfV{\mathbf V}
\def\bfW{\mathbf W}
\def\bfX{\mathbf X}
\def\bfY{\mathbf Y}
\def\bfZ{\mathbf Z}

\def\frD{\mathfrak D}
\def\frL{\mathfrak L}
\def\frS{\mathfrak S}

\def\frg{\mathfrak g}
\def\frz{\mathfrak z}
\def\frm{\mathfrak m}

\def\bfw{\mathbf w}

\def\R {{\mathbb R }}
 \def\C {{\mathbb C }}
  \def\Z{{\mathbb Z}}
  \def\H{{\mathbb H}}
\def\K{{\mathbb K}}
\def\N{{\mathbb N}}
\def\Q{{\mathbb Q}}
\def\A{{\mathbb A}}

\def\T{\mathbb T}
\def\P{\mathbb P}

\def\G{\mathbb G}

\def\cD{\EuScript D}
\def\cL{\EuScript L}
\def\cK{\EuScript K}
\def\cM{\EuScript M}
\def\cN{\EuScript N}
\def\cR{\EuScript R}
\def\cW{\EuScript W}
\def\cY{\EuScript Y}
\def\cF{\EuScript F}
\def\cT{\EuScript T}
\def\cB{\EuScript B}
\def\cE{\EuScript E}
\def\cO{\EuScript O}
\def\cP{\EuScript P}
\def\cH{\EuScript H}

\def\bbA{\mathbb A}
\def\bbB{\mathbb B}
\def\bbD{\mathbb D}
\def\bbE{\mathbb E}
\def\bbF{\mathbb F}
\def\bbG{\mathbb G}
\def\bbI{\mathbb I}
\def\bbJ{\mathbb J}
\def\bbL{\mathbb L}
\def\bbM{\mathbb M}
\def\bbN{\mathbb N}
\def\bbO{\mathbb O}
\def\bbP{\mathbb P}
\def\bbQ{\mathbb Q}
\def\bbS{\mathbb S}
\def\bbT{\mathbb T}
\def\bbU{\mathbb U}
\def\bbV{\mathbb V}
\def\bbW{\mathbb W}
\def\bbX{\mathbb X}
\def\bbY{\mathbb Y}

\def\kappa{\varkappa}
\def\epsilon{\varepsilon}
\def\phi{\varphi}
\def\le{\leqslant}
\def\ge{\geqslant}

\def\B{\mathrm B}

\def\la{\langle}
\def\ra{\rangle}

\def\F{{}_2F_1}
\def\FF{{}^{\vphantom{\C}}_2F_1^\C}

\newcommand{\dd}[1]{\,d\,{\overline{\overline{#1}}} }

\def\lambdA{{\boldsymbol{\lambda}}}
\def\alphA{{\boldsymbol{\alpha}}}
\def\betA{{\boldsymbol{\beta}}}
\def\mU{{\boldsymbol{\mu}}}
\def\PI{{\boldsymbol{\Pi}}}

\def\1{\boldsymbol{1}}
\def\2{\boldsymbol{2}}

\def\Th{\mathrm{T\!h}}

\begin{center}
\bf\Large

On spherical unitary representations of groups of spheromorphisms of Bruhat--Tits trees

\bigskip

\large \sc Yury A. Neretin%
\footnote{The research was supported by the grants FWF, Projects P25142, P28421, P31591.}
\end{center}

{\small Consider an infinite homogeneous tree $\cT_n$ of valence
	$n+1$, its group $\Aut(\cT_n)$ of automorphisms, and the group
	$\Hie(\cT_n)$ of its spheromorphisms (hierarchomorphisms), i.~e.,
	the group of homeomorphisms of the boundary of $\cT_n$ that locally coincide
	with transformations defined by automorphisms. We show that
	the subgroup $\Aut(\cT_n)$ is spherical in $\Hie(\cT_n)$, i.~e.,
	any irreducible unitary representation of $\Hie(\cT_n)$
	contains at most one $\Aut(\cT_n)$-fixed vector. We  present a combinatorial
	description of the space of double cosets of $\Hie(\cT_n)$ with respect to
	 $\Aut(\cT_n)$
	 and construct a 'new' family of spherical representations of $\Hie(\cT_n)$.
 We also show that the Thompson group $\Th$ has $\PSL(2,\Z)$-spherical unitary
representations.}

\section{Introduction}


\COUNTERS

\begin{figure}
	$$\epsfbox{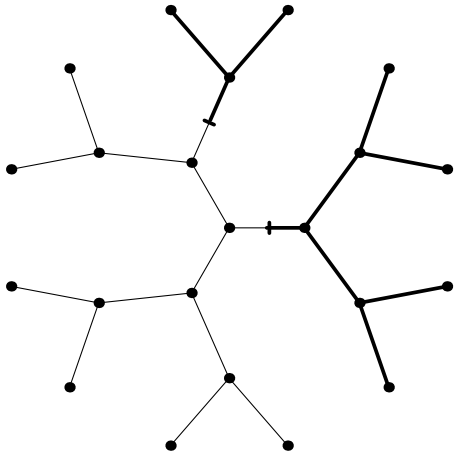}$$
	\caption{Ref. to Subsect. \ref{ss:1}. A piece of the Bruhat--Tits tree $\cT_2$. Transposing the thick branches
		we get an spheromorphism. }
	\label{fig:tree}
		\end{figure}

{\bf \punct Groups of spheromorphisms  of trees.%
\label{ss:1}} Fix an integer $n\ge 2$. The
{\it Bruhat--Tits tree} $\cT_n$ is the infinite tree such that each vertex belongs to $n+1$ edges,
see Fig. \ref{fig:tree}.	
Denote by $\Aut(\cT_n)$ the group of all automorphisms of $\cT_n$. It is a totally disconnected  locally compact group,
its topology is defined from the condition: stabilizers of finite subtrees are  open
in $\Aut(\cT_n)$.

Recall that Bruhat and Tits in 1966-1967 (see \cite{BT}) invented simplicial complexes
(Bruhat--Tits buildings), which are $p$-adic counterparts of noncompact Riemannian symmetric spaces.
Analogs of rank one noncompact symmetric  spaces (as the Lobachevsky plane) are Bruhat--Tits trees with $n$ being  powers of  prime $p$. In particular,
$p$-adic $\PSL(2)$ acts on the tree $\cT_p$. This fact became an initial point for investigations of group acting 
on trees, see, e.g., Tits \cite{Tits}, Serre \cite{Ser2}. Cartier \cite{Car} observed that
the groups $\Aut(\cT_n)$
are interesting objects from the point of view of representation theory and non-commutative harmonic analysis,
and these groups are relatives of $\SL(2)$ over real and $p$-adic fields. G.~Olshanski established that $\Aut(\cT_n)$ are type $I$
groups \cite{Olsh-trees-1} and obtained a pleasant classification \cite{Olsh-trees-2} of irreducible unitary representations of $\Aut(\cT_n)$
(see an exposition in \cite{Figa}, see also the work \cite{Cat} on tensor products).

The {\it boundary} $\Abs(\cT_n)$ of $\cT_n$ is a totally disconnected compact set, for a prime $n=p$
it can be identified with a $p$-adic projective line. The group $\Aut(\cT_n)$ acts by homeomorphisms
of the boundary.
A {\it spheromorphism} (or {\it hierarchomorphism}) of $\cT_n$ is a homeomorphism $q$
of  $\Abs(\cT_n)$ such that for each point $y\in \Abs(\cT_n)$ there is its neighborhood $\cN(y)$, 
in which $q$ coincides with some $r_y\in \Aut(\cT_p)$. In other words, we cut a
finite number of mid-edges of the tree  and get a collection of finite pieces $W_i$ and infinite pieces $U_j$.
We forget finite pieces and choose embeddings $\theta_j:U_j\to \cT_n$ such that images $\theta_j$ are mutually disjoint
and cover the whole tree (may be) without a finite piece,
see  Fig. \ref{fig:tree}.
The group $\Hie(\cT_n)$ of all spheromorphisms
of the tree $\cT_n$ is a locally compact topological group (see, \cite{GL}). 
The topology is defined by the 
condition: the subgroup $\Aut(\cT_n)$ is open and closed ({\it clopen}) in $\Hie(\cT_n)$.
The (countable) space of cosets $\Hie(\cT_n)/\Aut(\cT_n)$ has a discrete topology%
\footnote{So we have a group $G=\Hie(\cT_n)$ and a subgroup $K=\Aut(\cT_n)$ such that $K$ is a continuous totally
disconnected group and the homogeneous space $G/K$ is discrete. Topologies of this kind arise in representation theory
of infinite symmetric groups, see \cite{Ner-symm}, Subsect. 3.7; a group with such a topology is used below
in Sect. \ref{s:ness}.}.

There is a well-known discrete group $\Th$ consisting
of spheromorphisms%
\footnote{We can imagine the Bruhat--Tits tree as drawn on the plane $\R^2$. Then we get a structure of a cyclically ordered set
on the boundary $\Abs(\cT_n)$. The Thompson group $\Th$ is the group of all  spheromorphisms preserving the cyclic order
on $\Abs(\cT_2)$.} defined by 1965 R. Thompson in 1965. Initially it was proposed as a counterexample, and it 
really has strange properties but
also it is an interesting positive object (see, e.g., \cite{GS}, \cite{CFP}, \cite{Imb}, \cite{Pen}, \cite{Kap-Ser}, \cite{Bur},
\cite{Fos}, \cite{FN}).

\sm

The groups $\Hie(\cT_p)$ were introduced in 1984, \cite{Ner1}--\cite{Ner2} with the following 
reasoning:

\sm

1) For prime $n=p$ the group  $\Hie(\cT_p)$ contains the group of locally analytic diffeomorphisms
of the $p$-adic projective line. 

\sm

2) Unitary representations of the group of diffeomorphisms of the circle partially survive for 
the groups $\Hie(\cT_n)$.

\sm

3) The groups $\Hie(\cT_p)$ have several families of unitary representations that are spherical (see below) with respect to
(noncompact) subgroup $\Aut(\cT_n)$; in Addendum we explain why this  property seems to be distinguished.

\sm

The topic of the present paper are unitary representations, we list
  some  references on a wider context.  The groups  $\Hie(\cT_p)$ are simple as  abstract
groups (Kapoudjian \cite{Kap}), uniformly simple (Gal, Gismatullin, \cite{GG}),
compactly generated (Caprace, De~Medts \cite{CdM})
compactly presentable (Le~Boudec \cite{LeB}),
they have nontrivial $\Z_2$-central extensions%
\footnote{It is interesting to find  unitary faithful unitary representations
of this extension.} constructed by Kapoudjian \cite{Kap-vir}.
They have no property (T) (Navas, \cite{Nav}%
\footnote{Notice that families of spherical representations of $\Hie(\cT_n)$
in the boson and fermion Fock spaces
constructed in \cite{Ner2} approximate the trivial one-dimensional representation}).
These groups
are simple locally compact groups that do not admit a lattice (Bader, Caprace, Gelander, and Mozes, \cite{BCGM}, this is the first example of such kind).
See Kapoudjian \cite{Kap1}, Sauer, Thumann, \cite{ST} on action of $\Hie(\cT_n)$ on CW-complexes. These groups can be 
included to families of relatives \cite{Ner3}, \cite{Led},  \cite{ST}.
It seems to the author that these groups being locally compact have 
various properties of infinite-dimensional (or 'large') groups%
\footnote{For instance, constructions of spherical representations
both in \cite{Ner1}, \cite{Ner2} and below in Section 4 are distinctive
construction for infinite-dimensional groups.
On the other hand, a parallel with infinite-dimensional groups 
also is incomplete, apparently the  groups $\Hie(\cT_n)$ have no trains in the sense
of \cite{Ner-book}.}.

\sm

{\bf\punct Sphericity.}
Let $G$ be a topological group, $K$ is its subgroup. Let $\rho$ be an irreducible unitary
representation of the group $G$ in a Hilbert space $H$. We say that a {\it representation $\rho$ is $K$-spherical} if
$H$ contains a unique upto a scalar factor nonzero $K$-fixed vector $v$
(the {\it spherical vector}). Its matrix element
$$
\Phi(g)=\la \rho(g)v,v\ra_H,\qquad \text{where $\|v\|^2=1$,}
$$
is called a {\it spherical function}. This function is automatically 
$K$-biinvariant, i.e.,
$$
\Phi(k_1gk_2)=\Phi(g)\qquad \text{for $g\in G$, $h_1$, $h_2\in K$.}
$$
In other words, a spherical function is defined on the double coset space
$K\setminus G/K$.

\begin{definition}
	\label{def} Let $G$ be a topological group, $K$
	a closed subgroup.
The subgroup $K$ is {\rm spherical} if 

\sm

{\rm A)} For any irreducible unitary representation of $G$ the subspace of $K$-fixed vectors has dimension
$\le 1$.

\sm

{\rm B)} There is a faithful unitary representation of $G$ and a vector $v$ such that the stabilizer of $v$ is $K$.
\end{definition}

\sm

{\sc Remark.} The second condition is necessary for the following reason.
 Quite often (if $K$ is not compact or 'heavy'
in the sense of \cite{Ner-book}) a restriction of any nontrivial
irreducible unitary representation of $G$ to $K$ has not $K$-fixed vectors at all.
More generally, if a vector $v$ is fixed by $K$, then quite often
$v$ is automatically is fixed by a certain larger group $\wt K\supset K$.
Such phenomena were widely used in classical ergodic theory after 
Gelfand, Fomin \cite{GF} and Mautner  \cite{Mau}. 
A detailed investigation of such phenomena for Lie groups were
done by Moore \cite{Moo} and Wang \cite{Wang1}, for $p$-adic groups
by Wang \cite{Wang1}--\cite{Wang}. Kaniuth, Lau \cite{Kan} and Losert \cite{Los} discussed stabilizers 
of  vectors in unitary representations of general
locally compact groups%
\footnote{In their terminology subgroups that can be stabilizers of vectors
	 'satisfy separation property'.}.
 \hfill $\square$

\sm

{\bf\punct The purposes of the paper.}
We prove the following statements.

\begin{theorem}
\label{th:1}
The subgroup $\Aut(\cT_n)$ is spherical in $\Hie(\cT_n)$.
\end{theorem}



\begin{proposition}
\label{pr:1}
 Let $\Phi_1(g)$, $\Phi_2(g)$ be $\Aut(\cT_n)$-spherical functions on $\Hie(\cT_n)$. Then
$\Phi_1(g)\,\Phi_2(g)$ is a spherical function.
\end{proposition}

For known spherical pairs $G\supset K$ (finite-dimensional and infinite-dimensional)
double coset spaces $K\setminus G/K$ admit
explicit descriptions.  
In Section 3, we  present such a description for the double coset space
$$
\Aut(\cT_n)\setminus \Hie(\cT_n)/\Aut(\cT_n).
$$
Double cosets correspond to $(n+1)$-valent graphs $\Gamma$ consisting 
of two disjoint trees $T_+$ and $T_-$ and a collection of edges 
connecting vertices of $T_+$ with vertices of $T_-$ (cf. 'tree pairs diagrams'
in \cite{Bur}).

In Section 4
we apply Nessonov's construction \cite{Ness} of representations of
infinite symmetric group to obtain  a 'new' family of spherical representations
of $\Hie(\cT_n)$.

\sm

Addendum contains some comments on problem of sphericity for locally compact groups.
We also show that the Thompson group $\Th$ has $\PSL(2,\Z)$-spherical representations.

\sm

{\bf\punct Some questions.}
Theorem \ref{th:1} implies the following questions.

\sm

1) Is it possible to classify $\Aut(\cT_n)$-spherical functions
on $\Hie(\cT_n)$?

\sm

2) Is $\Hie(\cT_n)$ a type $I$ group?

\sm

3) Is it possible a harmonic analysis on the space $\Hie(\cT_n)/\Aut(\cT_n)$
in some sense%
\footnote{This is not a question about the decomposition of
$\ell^2$ on this space, see Addendum, Proposition
A.\ref{pr:el-2}.}.

\sm

4) Let $\rho$ be a spherical representation of $\Hie(\cT_n)$, let $P$
be the operator of orthogonal projection to $\Aut(\cT_n)$-fixed line. Consider 
the closure $\Gamma_\rho$ of $\rho(g)$, where $g$ ranges in $\Hie(\cT_n)$,
in the weak operator topology. Obviously (see Lemma \ref{l:HM})
the semigroup $\Gamma_\rho$ contains $P$, therefore 
$\Gamma_\rho$ contains operators of the form $\rho(g_1)P\rho(g_2)$
with $g_1$, $g_2\in \Hie(\cT_n)$. Does it contain something else?

\sm 

{\sc Remark.} 
The analog of the group of spheromorphisms for
	$n=\infty$  and its unitary representations
	are topics of a separate paper \cite{Ner-infinite}.
	\hfill $\square$

\section{Sphericity%
\label{s:spher}}

\COUNTERS

{\bf\punct Notation.%
\label{ss:termin-1}}
A {\it way} in the Bruhat--Tits tree is a sequence of vertices
$a_j$ such that $a_i$ and $a_{i+1}$ are adjacent
and $a_{i+2}\ne a_i$ for all $i$. We say that ways $a_i$ and $b_j$ are
{\it equivalent} if $a_i=b_{i+k}$ starting some $i$.
The {\it boundary}
(the notation: $\Abs(T_n)$)  of $\cT_n$ is the space of classes of equivalent ways.

\begin{figure}
	$$\epsfbox{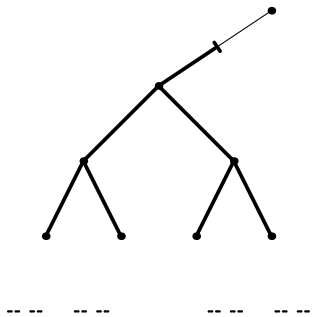}$$
	\caption{Ref. to Subsect. \ref{ss:termin-1}.
A branch of $\cT_2$ and the corresponding ball in $\Abs(\cT_2)$.}
\label{fig:ball}
\end{figure}

Let us cut the tree $\cT_n$ at the middle of an edge. We call two pieces
of the tree obtained in this way by {\it branches.}
Each branch $U$ determines a subset 
 $B=\Ba[U]$ in the boundary corresponding to ways lying in $U$.  
We call such subsets  by {\it balls}, see Fig. \ref{fig:ball}.
For a given ball $B$ denote by $\Br[B]$ the corresponding branch of the tree.
In particular, each mid-edge determines a partition
of $\Abs(\cT_n)$ into two disjoint balls. 
We define the {\it topology} on $\Abs(T_n)$ assuming
that balls  are clopen subsets in $\Abs(T_n)$,
this defines on $\Abs(\cT_n)$ a structure of totally disconnected
compact set. 

If $B_1$, $B_2$ are two balls, then
\begin{equation}
B_1\supset B_2,\quad\text{or}\quad B_2\supset B_1,\quad\text{or}\quad B_1\cap B_2=\varnothing
\label{eq:BB}
\end{equation}

\begin{lemma}
\label{l:union}
 Let $B_1\subset B_2\subset \dots$ be an increasing 
 sequence of balls. Then it has a maximal element or $\Abs(\cT_n)\setminus \cup_j B_j$
 is one point.
\end{lemma}

{\sc Proof.}
 Let a sequence of balls $B_j=\Ba[U_j]$ strictly decrease.
 Let $u_j$ be the corresponding mid-edges, and $[p_jq_j]$
 the corresponding edges, to definiteness assume $p_j\notin U_j$, $q_j\in U_j$.
 Then points $q_1$, $p_1$, $q_2$, $p_2$, \dots lye on a way. Let $a\in\Abs(\cT_n)$
be the limit of this way. Then $\cup B_j=\Abs(\cT_n)\setminus a$.
\hfill $\square$

\sm

We say that  $h\in \Aut(\cT_p)$ is {\it hyperbolic} if it has two fixed points $a$, $b$ on $\Abs(T_n)$
and induces a nontrivial shift on the two-side way $\dots x_{-1}$, $x_0$, $x_1$, \dots connecting $a$ and $b$.
Let $c$ be a point of the boundary. The {\it parabolic subgroup} $P_c\subset \Aut(\cT_n)$
is the group of transformations 
$g$ such that $g$ fixes $c$, and for any way $x_1$, $x_2$, \dots going to 
$c$ we have $gx_N=x_N$ for sufficiently large $N$.

\sm

{\bf \punct Proof of Theorem \ref{th:1}.}
The group $\Aut(\cT_n)$
has a normal subgroup $\Aut_+(\cT_n)$
of index 2 defined as follows.
Let us paint vertices of $\cT_n$  black and white
in such a way that each edge has edges of different colors.
The $\Aut_+(\cT_n)$ is
the subgroup of the group
preserving coloring. 
This defines a homomorphism 
of $\Aut(\cT_n)$ to the group
$\Z_2$ and therefore a one-dimensional representation of $\Aut(\cT_n)$.
Other nontrivial irreducible representations of $\Aut(\cT_n)$ are infinite-dimensional.
It is sufficient to prove the following
statement:

\begin{proposition}

 Consider an irreducible 
 unitary representation $\rho$ of $\Hie(\cT_p)$ in a Hilbert space $H$.
Denote by  $H^{\Aut_+}$  the subspace  of all $\Aut_+(\cT_n)$-fixed vectors.
 It is sufficient to prove that  
$\dim H^{\Aut_+}$ is $\le 1$.
\end{proposition}
 
 Denote by
$P$ the operator of orthogonal projection  to $H^{\Aut_+}$.
Clearly, 
\begin{equation}
P\rho(h)=\rho(h) P=P\quad \text{for all $h\in \Aut_+(\cT_n)$.}
\label{eq:Ph}
\end{equation}
For $g\in \Hie(\cT_n)$ we define an operator 
$\wh \rho(g):H^{\Aut_+}\to H^{\Aut_+}$ by
$$
\wh \rho(g):=P \rho(g) P.
$$
Clearly, $\wh\rho(g)$ depends only on a double coset
$ \Aut_+(\cT_n)\cdot g\cdot \Aut_+(\cT_n)$.

\begin{lemma}
 \label{l:P}
 The operators $\wt \rho(g)$ commute, i.e., for any
 $g_1$, $g_2\in \Hie(\cT_n)$ 
 \begin{equation}
 \label{eq:commutative}
 \wh \rho(g_1)\wh \rho(g_2)=\wh \rho(g_2)\wh \rho(g_1).
 \end{equation}
\end{lemma}

{\sc Reduction of Theorem \ref{th:1} to Lemma \ref{l:P}.}
Let the conclusion of the lemma hold. Assume that $\dim H^{\Aut_+}>1$.
Notice that $\wh\rho(g^{-1})=\wh\rho(g)^*$, therefore commuting bounded operators
$$
\wh\rho(g)+\wh\rho(g^{-1}), \qquad i\bigl(\wh\rho(g)-\wh\rho(g^{-1})\bigr),
$$
are self-adjoint. Therefore all operators $\wh \rho(g)$ have a proper  common
invariant subspace $V\subset H^{\Aut_+}$. Then $\Aut_+(\cT_n)$-cyclic span
of $V$ is a proper subspace in $H$. Indeed, let $v\in V$. Then
$$
P \rho(g)v=P \rho(g)Pv= \wh \rho(g)v\in V,
$$
and the projection of the cyclic span to $H^{\Aut_+}$ is contained to $V$.
\hfill $\square$

\begin{lemma}
	\label{l:HM}
 Let $h_j\in \Aut_+(\cT_n)$
 tend to infinity%
 \footnote{We say that $h_j$ {\it tends to $\infty$} if 
any compact subset of $\Aut_+(\cT_n)$ contains only a finite number of elements $q_j$.
In other words $h_j$ tends to infinity in the Alexandroff compactification of a locally compact space
$\Aut_+(\cT_n)$.}. Then for any unitary
representation $\rho$ of $\Aut_+(\cT_n)$ the sequence $\rho(h_j)$ converges 
to $P$ in the weak operator topology.
\end{lemma}

Equivalently for any nontrivial irreducible representation of $\Aut_+(\cT_n)$
the sequence $\rho(h_j)$ weakly converges to 0. For the group $\Aut(\cT_p)$
this holds for any irreducible unitary representation of dimension $>1$.
 This is proved in \cite{Lub}.
On the other hand this can be easily verified case-by-case starting Olshanski's classification theorem  \cite{Olsh-trees-2}.
Notice also that this is a counterpart of the well-known Howe--Moore theorem 
\cite{HM} about real Lie groups.

In fact, we need the following corollary.

\begin{corollary}
Let $h\in\Aut_+(\cT_n)$ be a hyperbolic element.
Then for any irreducible unitary representation $\rho$
of $\Aut_+(\cT_n)$ the sequence $\rho(h^m)$   
weakly converges to 0.
\end{corollary}

\sm

{\sc Proof of Lemma \ref{l:P}.} Fix a ball $B\subset\Abs(\cT_n)$.  Denote by
$G(B)$ the subgroup in $\Hie(\cT_n)$ consisting of spheromorphisms trivial outside $B$.
Clearly, 
 $$
   \Aut_+(T_n)\cdot G(B) \cdot \Aut_+(T_n)=\Hie(T_p),
   $$
 i.e., any double coset has a representative in $G(B)$.
 Choose two disjoint balls $B_1$, $B_2$.
 For a verification of (\ref{eq:commutative}) we can assume $g_1\in G(B_1)$, $g_2\in G(B_2)$.
 Choose a hyperbolic element $U\in \Aut_+(\cT_n)$ with an attractive fixed point $a\in B_2$.
 For $k>0$ we have
 $$
 U^k g_2 U^{-k}\in G(U^k B_2)\subset G(B_2).
 $$
 Hence $g_1$ and $U^k g_2 U^{-k}$ have disjoint supports, therefore they  commute. Thus,
 $$
 \rho(g_1) \, \rho(U^k)\, \rho(g_2)\, \rho(U^{-k})=  \rho(U^k)\, \rho(g_2)\, \rho(U^{-k})\, \rho(g_1).
 $$
 Multiplying this from the left and the right by $P$ and
 keeping in the mind (\ref{eq:Ph}), we get 
  $$
P \rho(g_1)  \rho(U^k) \rho(g_2) P=  P \rho(g_2) \rho(U^{-k}) \rho(g_1) P.
 $$
 Passing to the weak limit as $k\to+\infty$ and applying Lemma \ref{l:HM} we come to
   $$
P\, \rho(g_1)\,  P\, \rho(g_2)\, P=  P\, \rho(g_2)\, P\, \rho(g_1)\, P.
 $$
 This is the equality (\ref{eq:commutative}).
 \hfill $\square$
 
 \sm
  
{\bf\punct Proof of Proposition \ref{pr:1}.}

\begin{proposition}
\label{pr:2}
Let $G\supset K$ be topological groups. Assume that $K$ does not admit
nontrivial finite-dimensional unitary representations.
 Let $\Phi_1(g)$, $\Phi_2(g)$ be $K$-spherical functions on $G$. Then
$\Phi_1(g)\,\Phi_2(g)$ is a spherical function.
\end{proposition}

 \begin{lemma}
Let $\nu_1$ $\nu_2$ be unitary representations of a group $\Gamma$. If
the tensor product $\nu_1\otimes\nu_2$ contains a nonzero $\Gamma$-invariant vector,
then  both $\nu_1$ and $\nu_2$ have finite-dimensional subrepresentations.  
 \end{lemma}
 
{\sc Proof of the lemma.} Assume that an invariant vector exists.
 Denote the spaces of representations by $V_1$, $V_2$.
We identify $V_1\otimes V_2$ with the space of Hilbert--Schmidt operators
$V_1'\to V_2$, where $V_1'$ is the dual space to $V_1$. An invariant vector 
corresponds to an intertwining operator $T:V_1'\to V_2$.
The operator $TT^*$ is an intertwining operator in $V_2$. Since $TT^*$
is compact and nonzero, 
it has a finite-dimensional  eigenspace, and this subspace is $G$-invariant.
\hfill $\square$

 \sm
 
{\sc Proof of Proposition \ref{pr:2}.}
 Let $\rho_1$ and $\rho_2$ be $K$-spherical representations of $G$ in $H_1$ and $H_2$.
 Let $v_1$, $v_2$ be fixed vectors. By the lemma, $v_1\otimes v_2$ is a unique
 $K$-fixed vector in $H_1\otimes H_2$. The cyclic span
 $W$ of $v_1\otimes v_2$ is an irreducible subrepresentation.
 Indeed, let $W=W_1\oplus W_2$ be
 reducible. Then  both projections of $v_1\otimes v_2$ to $W_1$, $W_2$
 are $K$-fixed, therefore $v_1\otimes v_2$ must be contained
in one of summands, say $W_1$, and thus the cyclic span of $v_1\otimes v_2$
is contained in $W_1$, i.e., $W=W_1$.

Now we consider the representation of $G$ in $W$,
\begin{multline*}
\left\la\bigl(\rho_1(g)\otimes \rho_2(g)\bigr) v_1\otimes v_2,\, v_1\otimes v_2\right\ra_W
=\la \rho_1(g)v_1,v_1\ra_{H_1}\cdot\la \rho_2(g)v_2,v_2\ra_{H_2}
\\=
\Phi_1(g)\,\Phi_2(g).\qquad \square\!\!\!\!\!\!
\end{multline*}
  
{\sc Proof of Proposition \ref{pr:1}.}
Consider $\Aut(\cT_n)$-spherical 
representations $\rho_1$, $\rho_2$ of $\Hie(\cT_n)$. They
also are $\Aut_+(\cT_n)$-spherical.
Therefore their tensor product has 
a unique $\Aut_+(\cT_n)$-fixed vector.
This vector also is $\Aut(\cT_n)$-fixed.
\hfill $\square$

\section{The space of double cosets%
\label{s:coset}}

\COUNTERS

{\bf \punct Terminology.%
\label{ss:termin-2}}
Let $T$ be a tree, $A_1$, \dots, $A_N$
a collection of vertices. The {\it subtree spanned by}
$A_1$, \dots, $A_N$ is the minimal subtree containing these points.

Let $S$ be a finite tree. 
The {\it boundary} $\partial S$ of $S$
is the set of vertices of valence 1.

\begin{figure}
	$$\epsfbox{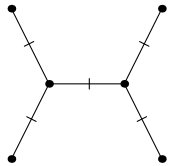}$$
	\caption{Ref. to Subsect. \ref{ss:termin-2}.
		The subdivision of the Bruhat--Tits tree.}
	\label{fig:subdivision}
\end{figure}

 We regard Bruhat-Tits trees as 1-dimensional complexes
with 0-cells located at vertices of the tree and mid-edges. Respectively, 1-cells are 
half-edges, see Fig. \ref{fig:subdivision}.

\sm

\begin{figure}
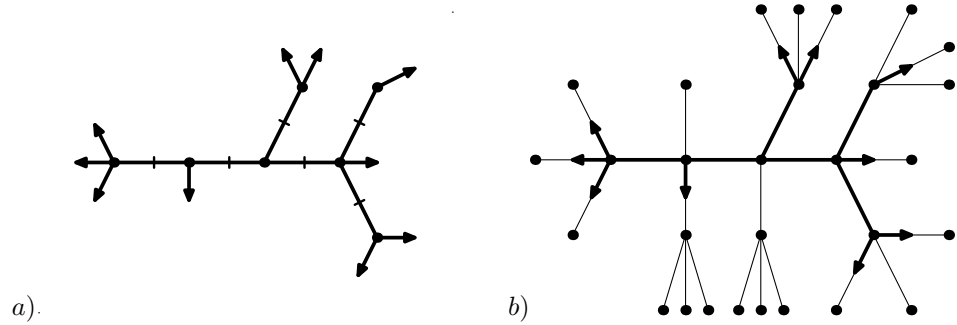

	$$a) \epsfbox{thorn.4}\qquad b)  \epsfbox{thorn.5}$$
	\caption{
		Ref. to Subsect \ref{ss:termin-2}.
		\newline
		a) A thorn (n=3). The left vertex is perfect.
		Cutting  the adjacent mid-edge off we get a reduced thorn. 
\newline	
b) A sub-thorn of the Bruhat-Tits tree $\cT_3$.
\newline
 On  Figure b) and figures  below we omit
mid-edges.
}
\label{fig:thorn}
\end{figure}

Let $R$ be a tree such that valences of all vertices
are $\le (n+1)$ and number of vertices is $\ge 3$.
A {\it thorn} $R$ is such a tree 
 equipped  with the following structure of an 1-dimensional simplicial 
complex.
Consider the subtree $R^\circ$ (the {\it skeleton of the thorn}) of $R$ spanned
by all vertices that are not contained in the boundary $\partial R$. Then 0-cells of the thorn
are vertices of $R$ and mid-edges of $R^\circ$. Respectively, 1-cells are half-edges of $R^\circ$ 
and edges of $R\setminus R^\circ$. We call vertices of $R^\circ$ by {\it vertices of thorn},
and points of  $\partial R$ by {\it spikes of the thorn}, see Fig. \ref{fig:thorn}.a.
\begin{figure}
	$$\epsfbox{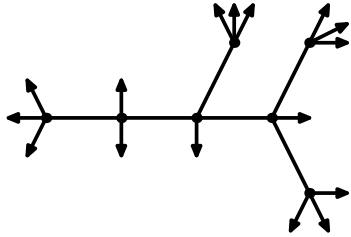}$$
	\caption{Ref. to Subsect \ref{ss:termin-2}. A perfect thorn ($n=3$).}%
	\label{fig:perfect}
\end{figure}
\begin{figure}
	$${\mathrm a)}\qquad\qquad\qquad\qquad\qquad\qquad
	{\mathrm b)}\qquad\qquad \epsfbox{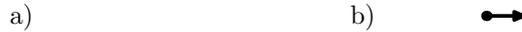}
	$$
	\caption{a) The empty thorn. b) The thorn with one vertex and one spike.}
	\label{fig:empty}
\end{figure}

Additionally, we allow  an {\it empty thorn} and a {\it thorn having 1 vertex and one spike},
see Fig. \ref{fig:empty}.

Denote by $\spike(R)$ the set of spikes of a thorn $R$, $\vert(R)$ the set
of vertices of  $R$.  

\sm

Two thorns $R_1$, $R_2$ are {\it isomorphic} if there is an
 isomorphism $R\to R'$ of complexes sending vertices to vertices
 and spikes to spikes.

Cutting a thorn in a mid-edge we get two {\it branches}.

\sm

We embed  thorns $R$ to the Bruhat--Tits tree $\cT_n$ 
  isomorphically sending vertices to vertices
  and spikes to mid-edges.
  We call  images of such  embeddings by  {\it sub-thorns}
  of the Bruhat--Tits tree, see \ref{fig:thorn}.b.

\sm

Let $R$ be a thorn.
We say a {\it thorn is  perfect} if all its vertices have valence
$(n+1)$, see Fig \ref{fig:perfect}. We say that a {\it vertex is perfect} if it is contained
in $\partial R^\circ$ and its valence is $(n+1)$, see \ref{fig:thorn}.b. More generally,
a {\it branch of a thorn is perfect} if all its vertices have valences $(n+1)$.

A thorn is {\it reduced} if it has no perfect vertices. Let $R$ be an arbitrary thorn.
Cutting of all perfect branches off we come to a reduced thorn (in particular, if $R$
is perfect, then the corresponding reduced thorn is empty.)

\sm 

{\bf \punct Clopen sets.}
Denote by $\Clop(\cT_n)$ the set of all nonempty clopen subsets
of $\Abs(\cT_n)$, by $\Clop^\circ(\cT_n)$ the subset
consisting of proper clopen subsets (i.e., we remove the point
of $\Clop(\cT_n)$
corresponding the whole $\Abs(\cT_n)$). 

  Clearly, any clopen subset $\Omega$ can be represented
as a union of a finite number of  disjoint balls
$$
\Omega:=B_1\sqcup\dots \sqcup B_\iota.
$$
 This representation is not unique,
since any ball $B$ can be canonically represented as a disjoint union of $n$ smaller balls.
It is easy to observe (see \cite{Ser2}, Addendum 'Structure of $p$-adic varieties',  or \cite{Ner2}), that the remainder
$\upsilon(\Omega)$
 of $\iota$ modulo $n-1$ is uniquely defined by $\Omega$.
 According this, $\Clop^\circ(\cT_n)$ splits
 as a disjoint union
 \begin{equation}
 \Clop^\circ(\cT_n)=\coprod_{\iota=0}^{n-2}  \Clop^\circ_\iota(\cT_n).
 \label{eq:iota}
 \end{equation}

\begin{proposition}
\label{pr:thorns}
	{\rm a)} Disjoint unions of balls $B_1\sqcup\dots \sqcup B_\iota$ are
	in one-to-one correspondence with sub-thorns of $\cT_n$.
	
		\sm 
	
	{\rm b)} Partitions $\Abs(\cT_n)=B_1\sqcup\dots\sqcup B_N$
	are in one-to-one correspondence with perfect sub-thorns of $\cT_n$. 
	
		\sm 
	
	{\rm c)} Nonempty clopen sets in $\Abs(\cT_n)$ are in one-to-one
	correspondence with reduced sub-thorns of $\cT_n$.
	
	\sm 
	
	{\rm d)} Orbits of $\Aut(\cT_n)$ on  $\Clop(\cT_n)$
	are numerated by equivalence classes of reduced thorns.
\end{proposition} 

\begin{figure}
 	\epsfbox{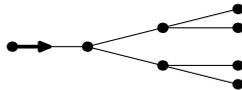}
 	\caption{A spike and the corresponding branch ($n=2$).}
 	\label{fig:spike}
 \end{figure}

{\sc Description  of the correspondence.}
Let $p$, $q$ be adjacent vertices of $\cT_n$.
Denote by $\overrightarrow{pq}$ the thorn having one vertex
$p$ and one spike in the mid-edge $pq$. Cutting the edge $pq$
at the mid-point we get two branches. We
 choose the branch $U$ containing $q$
 and the corresponding ball $B[\overrightarrow{pq}]$, see 
 Fig. \ref{fig:spike}.
 
 \sm 
 
{\it A sub-thorn $\longrightarrow$ a union of balls.} Consider  a sub-thorn in $\cT_n$.
 Then each spike corresponds to a ball.
 Taking a union of these balls we get a clopen subset with a given partition 
 into balls.
 
 Notice, that starting a perfect thorn we get the whole boundary $\Abs(\cT_n)$.

\sm

{\it A union of balls $\longrightarrow$ a sub-thorn.}
Conversely,
 fix a representation of $\Omega$ as a disjoint union
of balls $B_1\sqcup \dots \sqcup B_m$. Let
$U_1$, \dots $U_m$ be the corresponding branches of $\cT_n$.
Let $u_1$, \dots, $u_m$ be mid-edges that cut  these branches off.
We consider the minimal sub-thorn $R$ of $\cT_n$ containing 
$u_1$, \dots, $u_m$. 

\sm 

{\it A clopen set $\longrightarrow$ a  reduced sub-thorn.} Let $\Omega$ be a proper 
clopen set.
By Lemma \ref{l:union}, any sub-ball $B\subset \Omega$
is contained in a unique maximal sub-ball $\wt B\subset\Omega$.
We take the partition of $\Omega$ into maximal sub-balls
and take the corresponding thorn. Clearly,  it is reduced.
\hfill $\square$




\sm


{\bf \punct Double cosets and bi-thorns.} A {\it bi-thorn}
is the following collection of data $\{R,Q;\theta\}$:

\begin{figure}
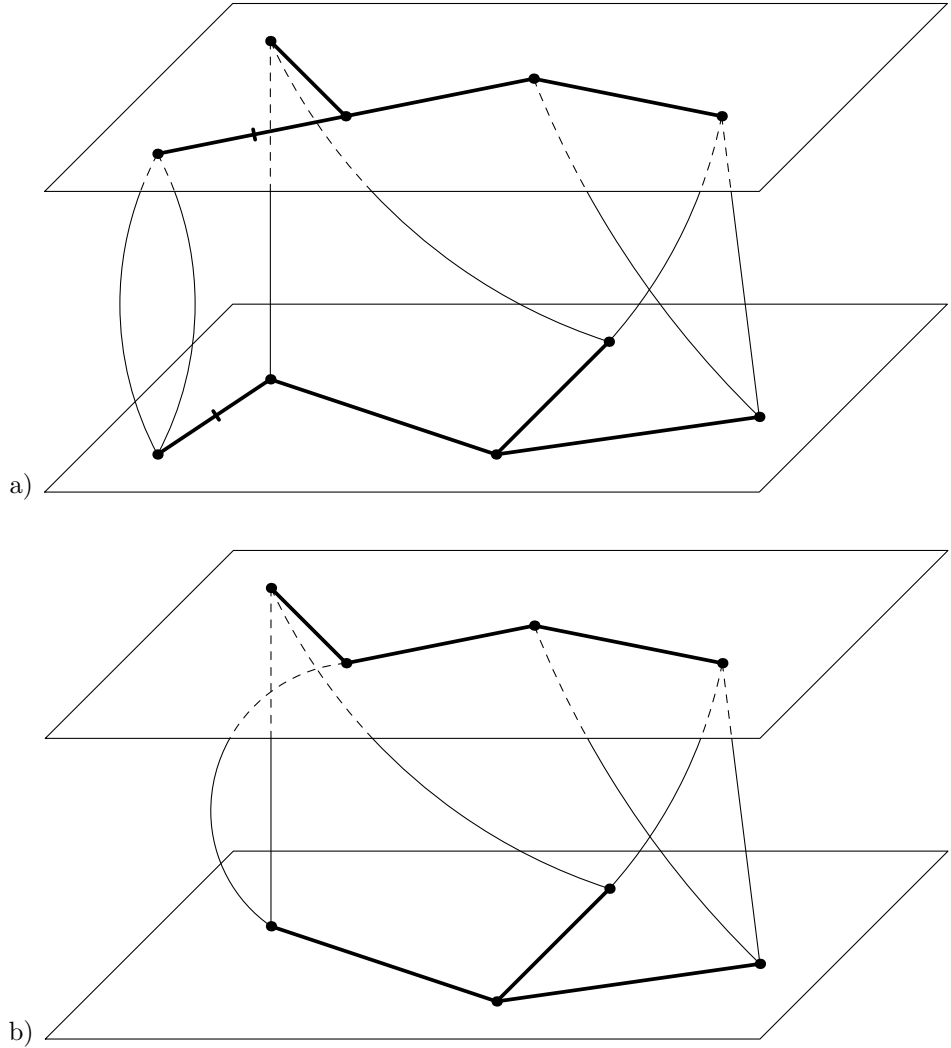

$$\!\! {\mathrm a)}\,\,	\epsfbox{thorn.7}$$

$$\!\! {\mathrm b)}\,\,	\epsfbox{thorn.8}$$
\caption{a) A bi-thorn. The left vertices of the upper and lower thorns are similar.
\newline
b) We cut off the left vertices and get a minimal bi-thorn
(an additional 'vertical' arc appears instead of two cut vertical arcs).
}
\label{fig:bi-thorn}
\end{figure}

\sm

$\bullet$ an ordered pair of perfect thorns $R$, $Q$ with the
same number of vertices;

\sm

$\bullet$  a bijection $\theta:\spike(R)\to \spike(Q)$.

\sm

We admit an {\it empty bi-thorn}. 

\sm

Equivalently, we have an $(n+1)$-valent graph $[R,Q;\theta]$,
which contains a pair of disjoint subtrees $R$, $Q$
and the remaining edges connect vertices of $P$ and vertices of $Q$
(we admit several edges between two vertices), see Fig. \ref{fig:bi-thorn}.

\sm

Consider a bi-thorn $\{R,Q;\theta\}$. 
Let $a$ be a vertex of $\partial (R^\circ)$,
$a'$ be a unique adjacent vertex of $R^\circ$.
Let
 $b$ a vertex
of $\partial (Q^\circ)$ and $b'$ the 
adjacent vertex. We say that $a$, $b$ are {\it similar }
if $\theta$ sends all spikes at $a$ to spikes at $b$, see Fig. \ref{fig:bi-thorn}.
In this situation, we can cut the mid-edges of $a'a$ and $b'b$.
The thorn splits into two pieces. We remove the piece
with two vertices $a$ and $b$ and modify $\theta$
saying that it sends the mid-edge of $a'a$ to 
the mid-edge of $b'b$. In this way we get a new thorn.

We say that a bi-thorn is {\it minimal} if it has not a pair of 
similar vertices.

\sm




\begin{proposition}
	\label{pr:bi-thorns}
	There is a canonical one-to-one correspondence between
	the double coset space $\Aut(\cT_n)\setminus \Hie(\cT_n) /\Aut(\cT_n)$
	and the set of minimal bi-thorns.
\end{proposition} 

Let us construct the correspondence.
Let $g\in\Hie(\cT_n)$. Take a ball $B=\Ba[U]$ and assume
that $gB$ is a ball,
 $gB=\Ba[V]$. We say that $g$ {\it regards} the ball $B$
 if the map $g:B\to gB$
 is induced by   an isomorphism
 of the branches $U\to V$.
 
 Let $g$ regard a ball $B$. Then there is a unique
 maximal ball $C=\wt B\supset B$ regarded by  $g$.
  Thus we get a partition
 $$\Abs(\cT_n) = C_1\sqcup C_2\sqcup \dots \sqcup C_N $$
 consisting of maximal balls regarded by
 $g$ and the corresponding partition  
$$\Abs(\cT_n) = gC_1\sqcup gC_2\sqcup \dots \sqcup gC_N $$
consisting of balls regarded by  $g^{-1}$.
We take thorns $R$ and $Q$ corresponding to this partitions,
by construction $g$ determines a bijection between their spikes.
\hfill $\square$

\begin{corollary}
\label{cor:finite}
 Fix $g\in\Hie(\cT_n)$. Fix an $\Aut(\cT_n)$-orbit $\cO$ 
 in $\Clop^\circ(\cT_n)$.
 Then for all but a finite number of elements  $\Omega\in \cO$ we have $g \Omega\in\cO$.
\end{corollary}

{\sc Proof.}
According the previous proof, $g$ canonically determines a pair of  sub-thorns $R$ and $Q$
of the Bruhat--Tits tree. The orbit $\cO$ corresponds to a certain reduced thorn
$T$. Elements $\Omega$  of the orbit correspond to sub-thorns $S$ in $\cT_n$ isomorphic to $T$.
Clearly, if $S\cap R=\varnothing$, then $g \Omega\in \cO$.
\hfill $\square$

\section{A family of spherical representations%
\label{s:ness}}

\COUNTERS

{\bf \punct The infinite symmetric group with Young subgroup.%
\label{ss:Ness}}
Fix $k$. Consider $k$ countable sets $\Pi_1$, \dots, $\Pi_k$ 
 and their disjoint union
$$
\PI:=\Pi_1\sqcup\dots\sqcup \Pi_k.
$$
First, consider
the group $G$ of all finitely supported permutations of $\PI$
and its (Young) subgroup $K$ preserving each $\Pi_j$.
Then $G\supset K$ is a spherical pair and according Nessonov
\cite{Ness} (see also, \cite{Ner-symm}, Sect. 8) all $K$-spherical functions on $G$
 have the following form $\Phi_S$. Consider a positive (semi)definite matrix $S$
 of size $k\times k$ with $s_{jj}=1$. Then
 $$
 \Phi_S(\sigma)=\prod_{p,q=1}^k
 s_{pq}^{\theta_{p,q}(\sigma)},\qquad \sigma\in G
 $$
 where $\theta_{p,q}(\sigma)$ is the number of elements $\alpha\in \Pi_p$ such that
 $\sigma\alpha\in \Pi_q$.
 
 To construct the corresponding unitary representations of $G$ we consider a Euclidean
 space $V$ and a collection of unit vectors $e_1$, \dots, $e_k$
 such that $\la e_p,e_q\ra_V=s_{p,q}$ (we can assume that $V$ is spanned by these vectors).
 Consider the tensor product%
 \footnote{Recall that a definition of a tensor product
 of an infinite family $H_j$ of Hilbert spaces requires a fixing
 of a distinguished unit vector $\xi_j\in H_j$ in each factor, a tensor product
 depends on a choice of $\xi_j$. For details, see, e.~g., \cite{Gui}, Appendix A.}
 $$
 \bigotimes_{p=1}^k\biggl( \bigotimes_{\alpha\in \Pi_p} (V, e_p)\biggr),
 $$
 we see that  factors are enumerated by elements of the set $\PI$.
The group $G$ acts by permutations of the factors. A unique $K$-fixed vector
is
$$
\cE:=\otimes_{p=1}^k e_p^{\otimes \infty}.
$$
The $G$-cyclic span of the vector $\cE$ is an irreducible spherical representation
of $G$.
 
Second, we notice that our representation can be extended by the continuity
to a larger group $\bfG$.  
 It consists of all permutations $\sigma$
of the set $\Pi$ such that for all $p$ for all but a finite number
of $\alpha\in \Pi_p$, we have $\sigma \alpha\in \Pi_p$
(permutations of factors in the tensor product are well-defined for such $\sigma$).
The spherical subgroup $\bfK$ consists of all permutations
preserving each subset $\Pi_p$.

\sm

{\bf \punct Embeddings of $\Hie(\cT_p)$ to the group $\bfG$.}
Consider a collection of reduced thorns $T_1$, \dots, $T_N$,
let they correspond to the same $\iota$ in the decomposition
(\ref{eq:iota}). Consider the corresponding $\Aut(\cT_n)$-orbits $\cO_1$, \dots, $\cO_N$
in $\Clop_\iota^\circ(\cT_n)$ and the complement $\cP$ to the union of these orbits.
Thus we get a partition
$$
\Clop_\iota^\circ(\cT_n)= \cP\sqcup \cO_1\sqcup\dots\sqcup \cO_N.
$$
Consider the group $\bfG$ corresponding to this partition.
By Corollary \ref{cor:finite}, the group $\Hie(\cT_n)$ is contained in $\bfG$. Obviously,
$\Aut(\cT_n)\subset \bfK$. So we can apply the Nessonov construction.

\sm

{\sc Remark.} Fix $\iota=0,$ 1, \dots, $n-2$.
Consider a Hilbert space $V$ and a countable set of unit vectors $e_T$
enumerated by reduced thorns whose number of spikes is $\iota$ modulo $n-1$.
Let this set have a unique limit point $e$ (and hence a sequence composed
of $e_S$ in any order converges to $e$. For a clopen subset $\Omega$ denote 
by $T(\Omega)$ the corresponding reduced thorn. Consider the following tensor product
$$
\cH:=\bigotimes_{\Omega\in \Clop_\iota^\circ} \bigl(V, e_{T(\Omega)}\bigr).
$$
The action of the group $\Hie(\cT_n)$ in $\cH$ by permutations of factors is well-defined
iff the following product absolutely converges for all hierarchomorphisms $g$:
$$
\Phi(g)=
\prod_{\Omega\in \Clop_\iota^\circ} \la e_{T(g\Omega)}, e_{T(\Omega)}\ra_V.
$$
Clearly, if the sequence $e_T$ converges fast enough, then this is the case.
In this situation, we get a spherical representation of $\Hie(\cT_n)$ in
$\cH$ with the spherical vector $\otimes_{\Omega\in \Clop_\iota^\circ} e_{T(\Omega)}$
and the spherical function $\Phi(g)$.

It can be interesting to find precise conditions for a family $e_T$ providing 
well-definiteness of this construction.

\section*{Addendum. Several comments on the sphericity phenomenon}


\sm

{\bf A.1. General remarks on sphericity.%
	\label{ss:general-spherical}}
Thus $\Aut(\cT_n)$ is a noncompact spherical subgroup in a locally compact group
$\Hie(\cT_p)$. According \cite{Ner-psl}, the subgroup $\PSL(2,\R)$
is spherical in the group $\Diff^3(S^1)$ of $C^3$-diffeomorphisms of the circle $S^1$.
We explain why this seems distinguished.

Phenomenon of sphericity was discovered by Gelfand in 1950, \cite{Gel}.
He showed that maximal compact subgroups $K$ in semisimple Lie groups
$G\supset K$ are spherical (as $\GL(n,\R)\supset \O(n)$ or
$\Sp(2n,\R)\supset \U(n)$).  Also symmetric subgroups in semisimple compact Lie groups
are spherical (as $\U(n)\supset \O(n)$ or $\O(2n)\supset\U(n)$).
Related spherical representations played a distinguished role in theory of unitary representations,
and spherical functions were an important standpoint for development of
modern theory of multi-dimensional special functions.

In 1979 Kr\"amer \cite{Kra} observed that simple compact Lie groups
have smaller spherical subgroups as $\O(2n+1)\supset \U(n)$ or $\Sp(2n+2)\supset \Sp(2n)\times \SO(2)$,
in the most of cases such pairs can be obtained from a Gelfand
pair $G\supset K$ by a minor enlargement of $G$ or minor reduction of $K$.
Mikityuk and Brion extended the Kr\"amer classification to semisimple compact groups.

There is also a story with finite spherical pairs $G\supset K$, see, e, g., \cite{Cec} 

\sm

On the other hand infinite-dimensional limits of Gelfand pairs (as $\GL(\infty,\R)\supset \O(\infty)$)
are spherical. G.~Olshanski \cite{Olsh-GB}, \cite{Olsh-symm} understood that such pairs have a substantial representation theory,
later there appeared related harmonic analysis. For infinite-dimensional (large) groups the
phenomenon of sphericity is more usual than for Lie group, and at least representation theory
can be developed in quite wide generality, see, e.~g. \cite{Ness0}, \cite{Ness}, \cite{Ner-symm}, \cite{Ner-spher-infty},
in Subsection \ref{ss:Ness} we used a construction of this kind.
In a known zoo, spherical subgroups are 'heavy groups' in the sense 
of \cite{Ner-book} (as the complete unitary group, the complete symmetric group,
the group of all measure preserving transformations).

\sm

Two examples mentioned in the beginning of the present subsection
are outside these two families. In one case a noncompact Lie group $\SL(2,\R)$
is a spherical subgroup in an infinite-dimensional group $\Diff^3(S^1)$,
in another case a noncompact subgroup $\Aut(\cT_n)$ is spherical in a locally compact group
$\Hie(\cT_n)$.

\sm

{\bf A.2. On compactness of stabilizers of vectors in unitary representations.%
	\label{ss:stabilizers-compact}}
In any case in substantial theory of unitary representations of Lie groups spherical subgroups
(in the sense formulated in Introduction) 
must be compact. There is a theorem of Moore \cite{Moo} about possible stabilizer of vectors 
in unitary representation, whose precise formulation is 
slightly sophisticated. We formulate a simpler statement.

Let $G$ be a connected Lie group, $Z$ the center; denote by $\frg\supset \frz$ their Lie algebras.
Denote by $\Ad_\frg(\cdot)$ the adjoint representation of $G$ in $\frg$, in fact we have a representation
of the quotient group $G/Z$ in the group $\GL[\frg]$ of all linear operators of the space $\frg$.

Let $\rho$ be a faithful irreducible unitary representation of $G$ in a Hilbert space $V$.
An irreducible faithful representation determines an injective homomorphism from $Z$ to the unit
circle $\T$ on the complex plane. For this reason $\dim\frz\le  1$, 
and we have 3 possibilities: $Z=\T$, $Z$ is finite, $Z$ is a dense subgroup in $\T$.

\begin{propositionA}
	\label{pr:moore}
	Let $G$, $\rho$, $V$ be as above.

	\sm
	
	{\rm (i)} Let  the image of $G/Z$ in the group $\GL[\frg]$  be closed.

	\sm
	
	{\rm (ii)} Let the center $Z$ be compact.
	
	\sm
	
	Then 
	
	{\rm a)} The stabilizer $K_v$ of a nonzero vector $v$ is compact.
	
	\sm
	
	{\rm b)} The stabilizer $L_v$ of the line $\C v$ is compact.

\end{propositionA}

{\sc Proof.} It is sufficient to prove the statement for the group $L_v$. By  definition
$L_v$ contains $Z$. Since $Z$ is compact, the image of $L_v$ in $G/Z$ is closed. Since
the $\Ad$-image of $G$ in $\GL[\frg]$ is closed, the $\Ad$-image of $L_v$ in  $\GL[\frg]$ also is closed.

We use Theorem 1.2 of Wang \cite{Wang1} (which is a strong version of the result of Moore \cite{Moo}).
We say that an element  $g\in \GL(\frg)$ is {\it pre-periodic} if it is semisimple%
\footnote{i. e., it is diagonalizable after a pass to the complexification.}
and its eigenvalues $\theta_j$
satisfy $|\theta_j|=1$. Equivalently, the closure of the set $\{g^m\}$, where $g$ ranges in $\Z$, is compact.
By \cite{Wang1}, for any  $g\in L_v$ there is a 
 subgroup $M_g$ such that:

\sm

1) $M_v\subset K_v$;

\sm

2) denote by $\frm_g$ the Lie algebra of $M_g$; then the image of $\Ad(g)$ in $\frg/\frm$
is pre-periodic.

\sm

However, {\it if a normal subgroup fixes a vector $v$, then it acts trivially on the whole space.} 
Indeed, let $r\in G$, $m\in M_g$. Then
$$
\rho(m)\, \rho(g)\, v=\rho(g)\, \rho(g^{-1} mg)\, v=\rho(g)v.
$$

Our representation is faithful and therefore the subgroup $M_g$ is trivial.
Thus the image $L_v/Z$ of $L_v$ in $\GL[\frg]$ is closed and consists of pre-periodic elements. It is more or less clear that
$L_v/Z$ is compact%
\footnote{To avoid a proof, we can refer to Lemma 1.3 from \cite{Wang} about a group with
	a {\it dense} set of
	pre-periodic elements.}. Since $Z$ is compact, $L_v$ also is compact.
\hfill $\square$

\sm

{\sc Remark 1.}
There are several reasons, for which we can not simply say: for unitary representations stabilizers of vectors (lines)
are compact. 

\sm 

a) Obviously we must consider faithful representations, since any closed normal subgroup $H\subset G$ can be a kernel of a representation.

\sm 

b) More serious sources of problems are twinnings.
Consider the group $\Isom(2)$
of orientation preserving isometries of the Euclidean plane.
Denote $Q:=\Isom(2)\times \Isom(2)$, we can regard an element of this group as a pair of matrices
of the form
\begin{equation}
\begin{pmatrix}
e^{it}&z\\
0&1
\end{pmatrix}, \qquad
\begin{pmatrix}
e^{is}&w\\
0&1
\end{pmatrix},\qquad\text{where $t$, $s\in \R$, $z$, $w\in\C$.}
\tag{A.1}
\end{equation}
Denote by $S\subset Q$ the subgroup consisting of pairs of matrices with $z=w=0$,
i.~e., $S=\SO(2)\times\SO(2)$. It is more-or-less obvious that $Q\supset S$ is a spherical pair
(the Wigner--Mackey trick, see, e.g., \cite{Kir}, 13.3, Theorem 1, immediately gives a classification
of irreducible unitary representations of $Q$). Next, choose an irrational real $\theta$ and take the subgroup
$G\subset Q$  consisting of pairs of matrices (A.1) satisfying the condition $s=\theta t$, consider the corresponding
subgroup $K=S\cap G$. The group $G$ is the {\it Mautner group} (see, e.~g. \cite{AM}).
Clearly, restricting an $S$-spherical representation of $Q$ to $G$ we get a $K$-spherical representation
of $G$. However, $K\simeq\R$ is not compact.

\sm

c) Consider the universal covering $G^\sim$ of the group $G=\SL(2,\R)$ and the universal covering
$R^\sim$ of the subgroup of rotations, $R^\sim\simeq\R$. Let $\rho$ be a faithful irreducible representation
of $G^\sim$ (see \cite{Puk}). Then $R^\sim$ has a discrete spectrum. For an eigenvector $v$ we have $L_v=R^\sim$
and $K_v\simeq \Z$. Both subgroups are non-compact. However, this non-compactness again
is artificial, in our case $L_v/Z$ is compact in $G^\sim/Z$.
\hfill $\boxtimes$

\sm

{\sc Remark 2.} If $G$ can be covered by a real algebraic group, then the conditions (i)-(ii) are fulfilled  
automatically. \hfill $\boxtimes$

\sm

Notice that for $p$-adic groups stabilizers of vectors in unitary representations
in interesting cases are compact (such stabilizers were topic of works of Wang \cite{Wang1}--\cite{Wang}).

\sm

{\bf A.2. The Mautner phenomenon for the groups for $\Hie(\cT_n)$.}
Let $\rho$ be a unitary representation of a group $G$. Assume that a subgroup $K$ fixes some vector $v$.
Then quite often $v$ is automatically fixed by certain larger group $\wt K$.
For $G=\Hie(\cT_n)$ we have the following statement:

\begin{propositionA}
\label{pr:mautner}
	Let $\rho$ be a unitary representation of $\Hie(\cT_n)$, let $v$ be a vector in the space
	of the representation.
	
	\sm
	
	{\rm a)} Let $h\in \Aut(\cT_n)$ be a hyperbolic element
	and $\rho(h)v=v$. Then $v$ is fixed by the whole subgroup $\Aut_+(\cT_n)$.
	
	\sm
	
	{\rm b)} Let $v$ be fixed by a parabolic subgroup $P_c\subset\Aut(\cT_n)$.
	Then $v$ is fixed by the whole subgroup $\Aut(\cT_n)$.
\end{propositionA}

This is obvious: nontrivial irreducible representations of $\Aut(\cT_n)$ have no fixed vectors
with respect to these subgroup (of course, this argument requires to look at Olshanski's
list \cite{Olsh-trees-2}). 
\hfill $\square$

\sm

{\bf A.3. A trivial spherical representation of $\Hie(\cT_n)$.}
Recall that the homogeneous space $\Hie(\cT_n)/\Aut(\cT_n)$ is countable and is equipped 
with the discrete topology. Therefore we have a quasi-regular representation
of $\Hie(\cT_n)$ in $\ell^2$ on this space, the natural orthogonal basis $\delta_z$ in $\ell_2$
is enumerated by points $z\in \Hie(\cT_n)/\Aut(\cT_n)$, the vector $\delta_z$ is the
delta-function supported by $z$.

\begin{propositionA}
\label{pr:el-2}
{\rm a)}
The representation of $\Hie(\cT_n)$ in $\ell^2\bigl(\Hie(\cT_n)/\Aut(\cT_n)\bigr)$
is irreducible and spherical, the spherical vector is $\delta_{z_0}$, where $z_0$ is the initial point of the homogeneous space,
the spherical function is {\rm 1} on $\Aut(\cT_n)$ and {\rm 0} outside this subgroup.

\sm

{\rm b)} Let $G$ be a topological group, $L$ a closed subgroup, let the 
homogeneous space $G/L$ be countable and discrete. Let $z_0$ be the initial point of $G/L$.
Let all orbits of $L$ on $G/L$ except $\{z_0\}$ be infinite. 
Then the representation
of $G$ in $\ell^2(G/L)$ is irreducible and spherical. The spherical vector is 
$\delta_{z_0}$ and the spherical function is zero outside $L$.
\end{propositionA}

{\sc Proof.} b) An $L$-invariant function on $G/L$ must be constant  
on orbits of $L$. Since a vector in $\ell^2$ can not have
infinite number of nonzero equal coordinates, we get that $\delta_{z_0}$
is the unique $L$-invariant vector. By the same argument as in the proof of Proposition
\ref{pr:2}, the $G$-cyclic span of $\delta_{z_0}$ is an irreducible subrepresentation
in $\ell^2$. However, this cyclic span contains all basis vectors $\delta_z$.

\sm

a) Keeping in mind Proposition \ref{pr:bi-thorns}, for any
element of $\Hie(\cT_n)/\Aut(\cT_n)$ we can assign a bi-thorn
$\{R,Q;\theta\}$ and an embedding of the thorn $Q$ to $\cT_n$.
The group $\Aut(\cT_n)$ acts preserving the bi-thorn and changing embeddings.
Clearly, if the bi-thorn $\{R,Q;\theta\}$ is non-empty, then orbits are infinite.
So, we can apply the statement b).
\hfill $\square$

\sm

{\bf A.4. A question about unitary representations of discrete groups.} 
It is well-know that questions about unitary representations of discrete groups quite often
are dangerous. By the Thoma theorem \cite{Tho2}, discrete groups are not type $I$ 
except groups that have  Abelian normal subgroups of finite index.  Absence
of type $I$ property implies numerous unpleasant phenomena (see, at least, the Glimm 
theorem \cite{Gli} about a bad Borel structure on the dual space). However, we formulate the following informal question.

\begin{questionA}
Consider a pair of countable groups $\Gamma\supset \Delta$ and let all orbits of $\Delta$
on $\Gamma/\Delta$ be infinite {\rm(}except the initial point{\rm)}. To find such pairs with
'interesting' $\Delta$-spherical representations of $\Gamma$. 
\end{questionA}

Apparently, interesting situations are rare.
However, there is a famous example of such a pair, which was basically discovered by in 1964 by Thoma \cite{Tho1} (see 
\cite{Olsh-symm}). We take the group $\Ss(\infty)$ of finitely supported permutations
of $\N$, let $\Gamma$ be $\Ss(\infty)\times \Ss(\infty)$ and $\Delta\simeq \Ss(\infty)$ be the diagonal subgroup. This was a start
of big story (representation theory of infinite symmetric groups), we only mention that in this case spherical representations can be extended
by continuity to a larger (continual) group (see \cite{Olsh-symm}, \cite{Ner-symm}).

The
pair of discrete groups $G\supset K$ from Subsect. \ref{ss:Ness} is spherical
(and again we have a continuous extension to a larger group $\bfG$).
A big zoo of examples of spherical representations in  \cite{Ner-symm} 
has a similar nature.

Next, consider the Thompson group $\Th$ realized as the group of all continuous 
piece-wise $\PSL(2,\Z)$-transformations of the real projective line $\R\P^1$, see \cite{Pen}, \cite{Imb}, by this
construction $\Th$ is embedded to $\Hie(\cT_2)$ and $\PSL(2,\Z)$ is contained in $\Aut(\cT_2)$.

\begin{propositionA}
Consider a unitary $\Aut(\cT_2)$-spherical representation $\rho$
of $\Hie(\cT_2)$ with  spherical vector $v$. Then the $\Th$-cyclic span
of $v$ is a $\PSL(2,\Z)$-spherical representation of $\Th$.
\end{propositionA}

{\sc Proof.} It sufficient to show that the restriction of $\rho$ to $\PSL(2,\Z)$
does not contain additional $\PSL(2,\Z)$-fixed vectors.
We take an hyperbolic element $h$ of $\PSL(2,\Z)$, say, $h=\begin{pmatrix}
                                                      2&1\\3&2
                                                     \end{pmatrix}$.
                   It is hyperbolic in $\Aut(\cT_2)$. By Proposition A.\ref{pr:mautner}.a,
 vectors fixed by $h$ are fixed by the whole group  $\Aut_+(\cT_2)$, and 
 a vector fixed by this subgrop is unique. 
 \hfill $\square$

 \tt
 \noindent
 Yury Neretin\\
Wolfgang  Pauli Institute/c.o. Math. Dept., University of Vienna \\
 \&Institute for Theoretical and Experimental Physics (Moscow); \\
 \&MechMath Dept., Moscow State University;\\
 \&Institute for Information Transmission Problems;\\
 yurii.neretin@math.univie.ac.at
 \\
 URL: http://mat.univie.ac.at/$\sim$neretin/

\end{document}